\documentclass[12
pt,twoside,a4paper]{amsart}
\usepackage{amssymb}
\usepackage[latin1]{inputenc}
\usepackage[T1]{fontenc}
\usepackage{times}

\def\hpq0{h^{p,q}_{\leq 0}}
\def\Hpq0{\H_{\leq 0}^{p,q}}

\def\dbar{\bar\partial}
\def\ddbar{\partial\dbar}

\def\R{{\mathbb R}}
\def\Rn{\R^n}
\def\C{{\mathbb C}}
\def\Cn{\C^n}

\def\Ker{{\rm Ker\,  }}

\def\H{{\mathcal H}}

\def\Re{{\rm Re\,  }}

\def\be{\begin{equation}}
\def\ee{\end{equation}}

\newtheorem{thm}{Theorem}[section]
\newtheorem{lma}[thm]{Lemma}

\theoremstyle{definition}

\theoremstyle{remark}

\newtheorem{preremark}{Remark}
\newtheorem{preex}{Example}

\numberwithin{equation}{section}

\begin{document}

\title[]
{ Subharmonicity properties of the Bergman kernel and some other functions
  associated to pseudoconvex domains.}

\author[]{ Bo Berndtsson}

\address{B Berndtsson :Department of Mathematics\\Chalmers University
  of Technology 
  and the University of G\"oteborg\\S-412 96 G\"OTEBORG\\SWEDEN,\\}

\email{ bob@math.chalmers.se}

\begin{abstract}
{ Let $D$ be a pseudoconvex domain in $\C^k_t\times\Cn_z$ and let
  $\phi$ be  a plurisubharmonic function in $D$. For each $t$ we
  consider the $n$-dimensional slice of $D$, $D_t=\{z; (t,z)\in D\}$,
let $\phi^t$ be the restriction of $\phi$ to $D_t$
  and denote by $K_t(z,\zeta)$ the Bergman kernel of $D_t$ with the
  weight function $\phi^t$. Generalizing a recent result of
  Maitani and Yamaguchi  (corresponding to $n=1$ and $\phi=0$)  we prove
  that $\log K_t(z,z)$ is a plurisubharmonic function in $D$. We also
  generalize an earlier results of  
  Yamaguchi concerning the Robin function and discuss similar results
  in the setting of $\Rn$. }
\end{abstract}

\bigskip

\maketitle

\section{Introduction}

 Let $D$ be a pseudoconvex domain in $\C_t^k\times\Cn_z$ and let
  $\phi$ be  a plurisubharmonic function in $D$. For each $t$ we
  consider the $n$-dimensional slice of $D$, $D_t=\{z; (t,z)\in
  D\}$ and the restriction,  $\phi^t$, of $\phi$ to $D_t$. Denote by $A^2_t= A^2(D_t, e^{-\phi^t})$ the Bergman space of
  holomorphic functions in $D_t$ satisfying
$$
\int_{D_t}|h|^2e^{-\phi^t} <\infty.
$$
The Bergman kernel $K_t(\zeta,z)$ of $A^2_t$ for a point $z$ in $D_t$ is
the unique 
holomorphic function of $\zeta$ satisfying
$$
\int_{D_t}h(\zeta)\overline {K_t(\zeta,z)}e^{-\phi(t,\zeta)}= h(z)
$$
for all functions $h$ in $A^2_t$. We shall prove the following
theorem.
\begin{thm} 
With the notation above, the function $\log K_t(z,z)$ is
plurisubharmonic, or identically equal to $-\infty$
in $D$. 
\end{thm}
In particular $\log K_t$ is plurisubharmonic in $t$ for $z$ fixed.
Theorem 1.1 was previously obtained in \cite{M-Y} in the case $n=1$
and $\phi=0$.

Theorem 1.1 may be seen as a complex version of Prekopa's theorem (
see \cite{Prekopa}) from convex analysis. This theorem says that if
$\phi(x,y)$ is a convex function in $\R^m_x\times\R^n_y$ and we define
the function $\tilde\phi$ in $\R^m_x$ by
\begin{equation}
e^{-\tilde\phi(x)}=\int_{\R^n} e^{-\phi(x,y)}dy,
\end{equation}
then $\tilde\phi$ is also convex. Equivalently, we may define
$$
\tilde\phi(x)=\log k(x),
$$
where
$$
k(x)=\left(\int_{\R^n} e^{-\phi(x,y)}dy\right)^{-1}.
$$
For each $x$ fixed, $k(x)$ can be seen as the ``Bergman kernel'' for
the space 
$
 \Ker(d)
$
of constant functions in $\R^n$, since the scalar
product in 
$$
L^2(\R^n,e^{-\phi(x,\cdot)})
$$
 of a function, $u$, with $k(x)$
equals the mean value of $u$, i e  the orthogonal projection of $u$ on
the space of constants. Thus Theorem 1.1 is what we get by replacing
the convexity hypothesis in Prekopa's theorem by plurisubharmonicity,
and the kernel of $d$ by the kernel of $\dbar$.  (In the complex
setting we also need to pay attention to the domains involved, since a
general pseudoconvex domain cannot be defined by an inequality
involving global plurisubharmonic functions.)

One interesting case of the theorem , where the analogy to Prekopa's theorem is more
evident,   is when $(t,0)$ lies in $D$ (for
$t$ in some open set), and $D_t$ and $\phi^t$ are both for fixed $t$
invariant under 
rotations $r_\theta(z)=e^{i\theta}z$. It then follows from the mean
value property for holomorphic functions that $K_t(\zeta,0)$ is
for each fixed $t$ a constant independent of $\zeta$,
$$
K_t=\left(\int_{D_t} e^{-\phi^t}\right)^{-1}.
$$
The following theorem from \cite{B1} is therefore a corollary of
Theorem 1.1.

\begin{thm}
Assume that for each fixed $t$, $D_t$ and $\phi^t$ are invariant under
rotations $r_\theta(z)=e^{i\theta}z$. Define the function $\tilde\phi$
by
$$
e^{-\tilde\phi(t)}= \int_{D_t} e^{-\phi(t,\xi)}.
$$
Then $\tilde\phi$ is plurisubharmonic.
\end{thm}

In particular, taking $\phi=0$ it follows that under the hypothesis of
Theorem 1.2, the function
$$
-\log |D_t|,
$$
where $|V|$ stands for the volume of a set, is plurisubharmonic. This
has recently been used by Cordero-Erausquin (see \cite{2Cordero}) to
give a proof of the Santal\'o inequality.

Still under the hypotheses of Theorem 1.2 we can also introduce a
large parameter, $p$, and define a function $\tilde\phi_p$ by
$$
e^{-p\tilde\phi_p(t)}= \int_{D_t} e^{-p\phi(t,\xi)}.
$$
Thus $e^{-\tilde\phi_p(x)}$ is the $L^p$-norm of
  $e^{-\phi(x,\cdot)}$. From the plurisubharmonicity of $\tilde\phi_p$ it is
  not hard to deduce that
$$
\tilde\phi_{\infty}=\inf_\xi\phi
$$
is also plurisubharmonic. 
This is one version of  Kiselman's {\it minimum principle} for plurisubharmonic
functions, \cite{Kis}. 

One main application of Kiselman's minimum principle, combined with a
use of the Legendre transform, was  to give a
procedure to ``attenuate the singularities'' of a given
plurisubharmonic function: Given an arbitrary plurisubharmonic
function $\phi$, and a number $c>0$,  Kiselman constructed a new
plurisubharmonic function 
which is finite at all points where the Lelong number of $\phi$ is
smaller that $c$ and still has a logarithmic singularity at points
where the Lelong number of $\phi$ exceeds $c$. This was in turn used  to
give an easy proof of Siu's theorem on the analyticity of sets defined
by Lelong numbers ( see  \cite{Kis3} ). 

It is a consequence of the Hörmander $L^2$-estimates for
the $\dbar$-equation that if $a$ is a point in a bounded domain
$\Omega$ and $\phi$ is plurisubharmonic in $\Omega$, then there is some
holomorphic function in $L^2(\Omega,e^{-\phi})$ which does not vanish
at $a$, if and only if the function $e^{-\phi}$ is locally integrable
in some neighbourhood of $a$. Using this we can prove the following
theorem, which can be seen as an alternative way of attenuating the
singularities of  plurisubharmonic functions. 

\begin{thm} Let $\Omega$ be a pseudoconvex domain in $\C^n$ and let
  $\phi$ be plurisubharmonic in $\Omega$. Let $\psi$ be the
  plurisubharmonic function in $\Omega\times\Omega$ defined by
$$
\psi(a,z)=\phi(z)+(n-1)\log |z-a|.
$$
Put
$$
\chi(a)=\log K_a(a,a),
$$
where $K_a$ is the Bergman kernel for $A^2(\Omega,e^{-2\psi^a})$. Then
$\chi$ is plurisubharmonic in $\Omega$, is finite at any point where
the Lelong number of  $\phi$ is smaller than 1 and has a logarithmic
singularity at any point where the Lelong number of $\phi$ is larger
than 1. The singularity set of $\chi$, $\{a; \chi(a)=-\infty\}$ is
equal to (the analytic) set where the Lelong number of $\phi$ is at
least 1.
\end{thm}

Theorem 1.3 suggests the introduction of a family of Lelong
numbers,
$$
\gamma_s(\phi,a)
$$
by
replacing the function $\psi$ by 
$$
\phi(z)+s\log|z-a|
$$
for $0\leq s<n$, and looking at points where the corresponding function
$\chi$ is singular. We would then get the so called integrability
index (see e g \cite{Kis2})  for $s=0$ and the classical Lelong
numbers for $s=n-1$.  

Theorem 1.1 is also intimately connected with another result
concerning curvature of vector bundles. We explain this in the
simplest case, when $D$ is the product $U\times\Omega$ of two domains
in $\C^k_t$ and $\C^n_z$ respectively. Let us also here assume that
$\phi$ is a bounded function, so that all the Bergman spaces
$A^2(\Omega,e^{-\phi^t})$ are equal as vector spaces, but the norm
  varies with $t$. We can then define a  vector bundle, $E$, over
  $U$ by taking $E_t=A^2(\Omega,e^{-\phi^t})$. This is then a trivial
    vector bundle, of infinite rank, with an hermitian metric defined
    by the Hilbert space norm. Our claim is that this vector bundle is
    positive in the sense of Nakano. This can be proved by methods
    very similar to the proof of Theorem 1.1. Such a  result however
    seems to be more natural in the setting of complex fibrations with
    compact fibers (so that the Bergman spaces are of finite
    dimension) and we will come back to it in a future publication.

We shall give two proofs of Theorem 1.1. The first, and simplest, one
is modeled on one proof of Prekopa's theorem given by Brascamp and
Lieb, \cite{Lieb}. Brascamp and Lieb used in their proof a 
version of Hörmander's $L^2$-estimates for the $d$-operator instead
of $\dbar$. They also proved directly this $L^2$-estimate by an
inductive procedure, using a version of Prekopa's theorem in smaller
dimensions. Our first proof adapts this proof to the complex case
but starts from Hörmander's theorem. 

The second proof does not use  Hörmander's theorem, but rather the a
priori estimates behind it. (It is somewhat similar to a recent proof
of Theorem 1.2 given by Cordero-Erausquin, \cite{Cordero}, which is in
turn inspired by \cite{BBN}.) Our proof is based on a representation of
the Bergman kernel as the pushforward of a subharmonic form. We have
included that proof since it 
seems to us that it will be useful in other similar situations. As an
example of that we give a generalization of a rather remarkable result of
Yamaguchi on the plurisubharmonicity of the Robin function,
\cite{Yam}. We finish the paper with a short discussion of what a
real variable version of a subharmonic form should be and how this
notion can be used to prove Prekopa's theorem and real variable
versions of Yamaguchi's result, \cite{Cardaliaguet}. 

I would like to thank Christer Borell for several interesting
discussions on the material of this paper.

\tableofcontents

\section{A special case of Theorem 1.1}

Let  $V$ be a smoothly bounded strictly pseudoconvex domain in $\Cn$, with
defining function $\rho$ so that $V=\{\zeta,\rho(\zeta)<0\}$. Let $U$
be a domain in 
$\C$ and let $\phi$ be a smooth strictly plurisubharmonic function in a
neighbourhood of $U\times V$. Fix a point $z$ in $V$ and let
$K_t(\cdot,z)$ be the Bergman kernel for $V$ with the weight function
$\phi^t$.The main step in
proving Theorem 1.1 is to prove  that in this situation, $K_t(z,z)$ is  a
subharmonic function of $t$.

For any square integrable holomorphic function $h$ in $V$ 
\begin{equation}
h(z)=\int_V h(\zeta)\overline{K_t(\zeta,z)}e^{-\phi^t}
\end{equation}
is independent of $t$. We shall differentiate this relation with
respect to $t$ and will then have use for the following lemma. 
\begin{lma}Let $V$ be a smoothly bounded strictly pseudoconvex domain in $\Cn$,
  and let $\phi$ be a function in $\Delta\times V$ which is smooth up
  to the boundary. Let $K_t(\zeta,z)$ be the Bergman kernel for the domain $V$
  with weight function $\phi^t$. Then $K_t$ is for $z$ fixed in $V$
  smooth up to the boundary of $V$ as a function of $\zeta$, and
  moreover depends smoothly on $t$.
\end{lma}
\begin{proof} Let $v_t$ be a smooth function in $V$ supported in a small
  neighbourhood of $z$, depending smoothly on $t$,  and put $f_t=\dbar
  v_t$. Let $\alpha_t$ be the 
  solution of the $\dbar$-Neumann problem
$$
\Box_t\alpha_t=(\dbar\dbar^*_t+\dbar^*_t\dbar)\alpha_t=f_t,
$$
where $\dbar^*_t$ is the adjoint of $\dbar$ with resepect to the
weight $\phi^t$. Since $u_t=\dbar^*_t\alpha_t$ is the minimal solution
in $L^2(V,e^{-\phi^t})$ to the equation $\dbar u=f_t$ we have
$$
u_t(\zeta) =v_t(\zeta)-\int_{\chi\in V} v_t(\xi) K_t(\zeta,\xi)e^{-\phi^t}.
$$
Choosing $v_t$ appropriately ( i e so that $v_t e^{-\phi^t}$ is a
radial function of integral one in a  small ball with center $z$) we
get that the last term on the right hand side is equal to $K_t(\zeta,z)$. 
It is therefore enough to prove that $\alpha$ has the smoothness
properties stated. To see this, note that if $t$ is close to $0$
$$
\Box_t=\Box_0- S_t,
$$
with $S_t$ an operator of order 1 with smooth coefficients which
vanishes for $t=0$. Hence
$$
(I-R_t)\alpha_t:=(I-\Box_0^{-1}S_t)\alpha_t=\Box_0^{-1} f_t.
$$
For $t$ sufficiently close to 0 we can invert the operator $I-R_t$ and
the lemma follows from basic regularity properties of the
$\dbar$-Neumann problem in strictly pseudoconvex domains.
\end{proof}

We now differentiate the relation 2.1 with respect to $\bar t$, using
the lemma. Let us denote by
$\partial^{\phi}_t$ the differential operator
$$
e^\phi\frac{\partial}{\partial t} e^{-\phi}=\frac{\partial}{\partial
  t}-\frac{\partial\phi}{\partial t}.
$$
It follows that the function
$$
u=\partial^{\phi}_t K_t
$$
is for fixed $t$ orthogonal to the space of holomorphic functions in
$A^2_t$. By the reproducing property of the Bergman kernel we have
$$
\Phi(t):=K_t(z,z)=\int_V
K_t(\zeta,z)\overline{K_t(\zeta,z)}e^{-\phi^t}.
$$
We shall use this formula to compute $
\partial^2 \Phi/\partial t\partial\bar t.$
We first get, using the notation $\dbar_t=\partial/\partial{\bar t}$ 
$$
\frac{\partial\Phi}{\partial\bar t}=\int_V\dbar_t
 K_t\overline{K_t}e^{-\phi^t} +  
\int_V K_t\overline{\partial^{\phi}_t K_t}e^{-\phi^t}.
$$
Since $K_t$ is holomorphic and $u$ is orthogonal to the space of
holomorphic functions, the second term vanishes. We next differentiate
once more. 
$$
\frac{\partial^2\Phi}{\partial t\partial\bar t}=
\int_V|\dbar_t K_t|^2e^{-\phi^t} +
\int_V \partial^{\phi}_t \dbar_t K_t
\overline{K_t}e^{-\phi^t} .
$$
Using the commutation rule  
\begin{equation}
\partial^{\phi}_t\dbar_t=
\dbar_t\partial^{\phi}_t
    +\phi_{t\bar t}
\end{equation}
in the second term we get 
$$
\frac{\partial^2}{\partial t\partial\bar t}\Phi=
\int_V|\dbar_t K_t|^2e^{-\phi^t} +\int_V
\phi_{t\bar t}|K_t|^2 e^{-\phi^t} +\int_V\dbar_t\partial^{\phi}_t
K_t\overline{ K_t}e^{-\phi^t} . 
$$
Moreover, by differentiating the relation 
$$
0=\int_V\partial^{\phi}_t K_t\,\overline{ K_t}e^{-\phi^t}
$$
we find that
$$
\int_V\dbar_t\partial^{\phi}_t\, K_t\overline{
 K_t}e^{-\phi}_t =-\int_V |\partial^{\phi}_tK_t|^2 e^{-\phi^t}=-\int_V
|u|^2 e^{-\phi^t}.
$$
All in all we therefore have that
\begin{equation}
 \frac{\partial^2\Phi}{\partial t\partial\bar t}=
\int_V|\dbar_t K_t|^2e^{-\phi^t} +\int_V
\phi_{t\bar t}|K_t|^2 e^{-\phi^t} -\int_V |u|^2e^{-\phi^t}.
\end{equation}

To estimate the last term we note that  $u$ solves the $\dbar$-equation
$$
\dbar u:= f=  \dbar\partial^{\phi}_t K_t =K_t
\dbar\frac{\partial\phi}{\partial t},
$$
(the last equation follows from a commutation rule similar to 2.2 since $K_t$
is holomorphic). Moreover, $u$ is the minimal solution to this
equation, since $u$ is orthogonal to the space of holomorphic
functions. By Hörmander's theorem ( see \cite{Demailly} for an appropriate formulation ) 
we therefore get that
$$
\int_V |u|^2 e^{-\phi^t}\leq \int_V \sum(\phi^t)^{j  \bar k}
f_j\bar f_k e^{-\phi^t}, 
$$
where $ (\phi^t)^{j  \bar k}$ is the inverse of the complex
Hessian of $\phi^t$. 
Inserting this into 2.3 and discarding the first (nonnegative) term we
have
$$
 \frac{\partial^2\Phi}{\partial t\partial\bar t}\geq
\int_V |K_t|^2  D  e^{-\phi^t},
$$
where
$$
D =\phi_{t\bar t} -\sum(\phi_{z_j\bar z_k})^{-1}\phi_{t \bar
  z_j}\overline{\phi_{t \bar z_k}}.
$$
$D$ equals precisely the determinant of the full complex
Hessian of $\phi$ divided by the determinant of the Hessian of
\(\phi^t\). Since $\phi$ is strictly plurisubharmonic, this quantity is positive,
and it follows that $\Phi$ is subharmonic.

To see that in fact even $\log K_t$ is subharmonic we change the
weight function $\phi$ to $\phi(t,\zeta) + \psi(t)$ where $\psi$ is an
arbitrary smooth subharmonic function. The Bergman kernel for the new weight
$\phi+\psi$ is $e^{\psi} K_t$, where $K_t$ is the Bergman kernel for $\phi$.
Therefore $e^{\psi} K_t$ is subharmonic for any choice of subharmonic
function $\psi$. This  implies that $\log K_t$ is subharmonic. 

\section{The general case of Theorem 1.1}
In the previous section we have proved Theorem 1.1 when the domains
$D_t$ are smoothly bounded and do not depend on $t$, under the extra
assumption that $\phi$ is smooth up to the boundary. The general case
is in principle a rather straightforward consequence of this special
case. There is however one subtility, arising from the
fact that some of the fiber domains $D_t$ may not be smoothly
bounded. This happens at points where  the topology of the fiber
changes, something which is not at all excluded by our hypotheses. (
The simplest such example is when  $D_t=\{\psi(z)<\Re t\}$ where
$\psi$ is a subharmonic function of one variable with two logarithmic
poles.When $\Re t$ is large negative, $D_t$ is a union of two disjoint
islands around the poles. The two islands come closer as $\Re t$
increases and eventually touch in  a figure eight, after which
they join to one single domain.)

\begin{lma}
Let $\Omega_0$ and $\Omega_1$ be bounded domains in $\Cn$, with
$\Omega_0$ compactly included in $\Omega_1$. Let $\phi_j$ be a sequence
of continuous weight functions in $\Omega_1$ such that
$$
\phi_j=\phi
$$
in $\bar\Omega_0$ and that $\phi_j$ increases and tends to to infinity
almost everywhere in 
$\Omega_1\setminus \Omega_0$. 
 Assume that the space of
holomorphic functions in $L^2(\Omega_1,e^{-\phi_0})$ is dense in the space of
holomorphic functions in $L^2(\Omega_0,e^{-\phi_0})$.
Fix a point $z$ in $\Omega_0$ and
let $K_j$ be the Bergman kernel for $z$ in $L^2(\Omega_1,\phi_j)$. Let
$K$ be the Bergman kernel for $z$ in  $L^2(\Omega_0,\phi)$.

Then $K_j(z,z)$ increases to $K(z,z)$.
\end{lma}

\begin{proof} The extremal characterisation of Bergman kernels,
$$
K(z,z)=\sup |h(z)|^2,
$$
where the supremum is taken over all holomorphic functions of
$L^2$-norm at most 1 makes it clear that $K_j(z,z)$ is an
increasing sequence and that each $K_j(z,z)$ is smaller than
$K(z,z)$. Since
$$
K_j(z,z)=\int_{\Omega_1}|K_j|^2e^{-\phi_j}
$$
it follows in particular that $K_j$ has uniformly bounded norm in
$L^2(\Omega_1, e^{-\phi_j}) $. The sequence $K_j$ therefore has a weakly
convergent subsequence in $L^2(\Omega_0, e^{-\phi}) $. Let $k$ be the
limit of some weakly convergent subsequence. If $h$ lies in
$L^2(\Omega_1, e^{-\phi_0} )$ we have that
$$
|\int_{\Omega_1\setminus\Omega_0}h\overline{K_j} e^{-\phi_j}|^2\leq
\int_{\Omega_1\setminus\Omega_0}|h|^2 e^{-\phi_j} \|K_j\|_{\phi_j}^2
$$
tends to zero. It follows that any weak limit $k$ satisfies
$$
h(z)= \int_{\Omega_0}h\bar k e^{-\phi}.
$$
Since holomorphic functions in $L^2(\Omega_1,e^{-\phi_0})$ are dense
 in $L^2(\Omega_0,e^{-\phi_0})$, the same relation holds for
any $h$ in  $L^2(\Omega_0,e^{-\phi_0})$. 
Since $k$ is necessarily also holomorphic, $k=K$ and the limit is in
fact uniform on compact subsets of $\Omega_0$. In particular
$$
\lim K_j(z)=K(z).
$$
\end{proof}
The proofs of the next two lemmas is similar but simpler and is
therefore omitted. 
\begin{lma} Let $\Omega$ be a bounded domain and $\phi$ a plurisubharmonic
  weight function . Let $\Omega_j$ be an increasing
  family of subdomains with union equal to $\Omega$. Let $z$ be a
  fixed point in $\Omega_0$ and let $K_j$ and $K$
  be the Bergman kernels for $\Omega_j$ and $\Omega$ (with weight
  function $\phi$) respectively. Then $K_j(z,z)$ decreases to
  $K(z,z)$.
\end{lma}
\begin{lma}
 Let $\Omega$ be a bounded domain and $\phi_j$ a decreasing sequence
  of  plurisubharmonic
  weight functions. Let $z$ be a
  fixed point in $\Omega$ and let $K_j$ and $K$
  be the Bergman kernels for the weight functions $\phi_j$ and $\phi$
  respectively. Then $K_j(z,z)$ decreases to $K(z,z)$.
\end{lma}

To verify one of the hypotheses in Lemma 4.1 we need an
approximation result.

\begin{lma} Let $\Omega_0$ and $\Omega_1$ be smoothly bounded
  pseudoconvex domains in $\Cn$ with $\Omega_0$ compactly included in
  $\Omega_1$. Assume there is a smooth  plurisubharmonic function $\rho$ in
  $\bar\Omega_1$ such that $\Omega_0=\{z\in\Omega_1,
  \rho(z)<0\}$. Then holomorphic functions in $L^2(\Omega_1)$ are
  dense in the space of holomorphic functions in $L^2(\Omega_0)$.
\end{lma}
\begin{proof}
Let $h$ be a square integrable holomorphic function in $\Omega_0$.
The crux of the proof  is to approximate $h$ by functions holomorphic in a
neighbourhood of the set $X=\{\rho\leq 0\}$. This can be done by
standard $L^2$-theory if 0 is a regular value of $\rho$ so that the
boundary of $\Omega_0$ is smooth. In the non-smooth case, the
possibility to approximate with function holomorphic near $X$ follows
from a result by Bruna and Burgues, cf Theorem B in \cite{Br-Bu}. 

Next, we let $h$ be holomorphic near $X$ and show
how to approximate $h$ with functions holomorphic in $\Omega_1$. Let
$H$ be an arbitrary extension of $h$ from a neighbourhood of $X$ to a
smooth function with compact 
support in $\Omega_1$ and put $f=\dbar H$. Let $k_j(s)$ be a sequence of
increasing convex functions that vanish for $s<0$ and tend to infinity
for $s>0$ and set 
$\phi_j=k_j\circ\rho.$ By Hörmander's theorem, \cite{Horm},  we can solve the
equation $\dbar v_j= f$ with estimates in
$L^2(\Omega_1,e^{-\phi_j})$. Since $f$ is supported in the complement
of $\Omega_0$ it follows that $v_j$ tends to zero in
$L^2(\Omega_0)$. Hence $H-v_j$ is an approximating sequence. 
\end{proof}

The final lemma gives the semicontinuity of $K_t$.
\begin{lma} Let $D=\{(t,z); \rho(t,z)<0\}$ where $\rho$ is smooth and
  strictly plurisubharmonic near the closure of $D$ and moreover has
  non-vanishing gradient on $\partial D$. Assume $\phi$ is smooth and
  plurisubharmonic near the closure of $D$.  Then $K_t(z,z)$ is for fixed
  $z$ upper semicontinuous as a function of $t$.
\end{lma}
\begin{proof}
Consider a point $t$ and let $s$ be nearby points tending to $t$. We
may choose $\epsilon>0$ so that all fibers $D_s$ are contained in the
open set $V$  where $\rho(t,z)<\epsilon$. 
Note that  the set-valued function $t\rightarrow D_t$
is lower semicontinuous, in the sense that if $D_t$ contains a compact
set $K$, the $K$ is contained in all $D_s$ for $s$ sufficiently close
to $t$. Let $K_s(\zeta,z)$ be the Bergman kernel of $D_s$ for a fixed
point $z$. Since   the domains $D_s$ all contain a fixed open
neighbourhood of $z$ the $L^2$-norms of $K_s$ are bounded. Any
sequence of $K_s$ therefore has a subsequence weakly convergent on any
compact subset of $D_t$. The $L^2$-norm of any weak limit $k$ can not
exceed the liminf of the $L^2$-norms of $K_s$ over $D_s$. By the
extremal characterization of Bergman kernels it follows that
$$
\limsup K_s(z,z)\leq K_t(z,z),
$$
so we are done.
\end{proof}

We can now complete the proof of Theorem 1.1, and start by proving
that $\log K_t$ is plurisubharmonic in $t$ for $z$ fixed. We first
assume that $D$ is smoothly bounded, defined as
$$
D=\{(t,z); \rho(t,z)<0\}
$$
where $\rho$ is smooth and strictly plurisubharmonic near the closure
of $D$. We also assume that $\phi$ is smooth and plurisubharmonic near
the closure of $D$. Assume first $k=1$ and 
fix a point $t$ in $\C$, say $t=0$. If
$U$ is a sufficiently small neighbourhood of 0 all
the fibers $D_t$ are contained in a fixed pseudoconvex domain
$V=\{\rho(0,\zeta)<\epsilon \}$. In $U\times V$ we can compose $\rho$
with an increasing sequence of smooth convex functions $k_j$ that tend
to infinity when $\rho$ is positive. We can now apply the result from
section 3 to $U\times V$ with $\phi$ replaced by $\phi_j=\phi +k_j\circ\rho$
and  let $j$ tend to infinity. Since the set where a smooth strictly
subharmonic function equals zero has zero measure, $\phi_j$ tends to
infinity a e in $\Omega_1\setminus\Omega_0$ By Lemma 4.1 it follows
that $\log K_t$ can be 
written as an increasing limit of functions subharmonic with respect
to $t$. Since, by the last lemma, $\log K_t$ is also upper
semicontinuous it follows that it is subharmonic. Again by the upper
semicontinuity we get that $\log K_t$ is plurisubharmonic if $k\geq 1$
since its restriction to any line is subharmonic. 

It is now easy to remove the extra hypothesis on $D$ and $\phi$. If
$D$ is an arbitrary pseudoconvex open set it has a smooth strictly
plurisubharmonic exhaustion function, and so can be written as an
increasing union of domains  of the type satisfying the extra hypotheses.
Near each such domain we can regularize $\phi$ by convolution. From
lemmas 4.2 and 4.3 we get that $\log K_t$ is a decreasing limit of
plurisubharmonic functions, and so is plurisubharmonic, or identically
equal to minus infinity.

We have thus proved that, under the hypotheses of Theorem 1.1, $\log
K_t$ is subharmonic as a function of $t$ for $z$ fixed. To see that it
is plurisubharmonic in $t$ and $z$ jointly we use, as in \cite{Yam},
the Oka trick of variation of the domain. We need to prove that, for
any choice of $a$ in $\Cn$, the function 
$$
\log K_t(z+ta,z+ta)
$$
is subharmonic in $t$. But, this is precisely the Bergman kernel at
$z$ for
the domain
$$
D_t- ta
$$
with the weight function translated similarily. Since the translated
domains are also pseudoconvex, and the translated weight function is
still plurisubharmonic, it follows that $\log K_t(z+ta,z+ta)$ is
subharmonic in $t$ and we are done.

\section{ Subharmonic currents}

We shall next give an alternate proof of Theorem 1.1 which is based on
a representation of the Bergman kernel as the pushforward of a
subharmonic form. To prepare for this we give in this section some
general facts on subharmonic forms or currents.

Let $T$ be a current of bidimension $(1,1)$, i e of bidegree $(n,n)$
in $U\times\C^n$ where $U$ is an open set in $\C$. We say that $T$ is
subharmonic if 

$$
i\ddbar T\geq0.
$$
Let $\pi$ be the projection from $\C_t\times\Cn_z$ to $\C_t$. 
If $T$ is compactly supported in the fiber direction, so that the
support of $T$ is included in $U\times K$ with $K$ a compact subset of
$\C^n$ the pushforward $\pi_{*}(T)$ of $T$ to $U$ is the distribution
in $U$ defined by 
$$
\pi_{*}(T).\chi= T.\pi^*\chi
$$
for any smooth compactly supported $(1,1)$ form $\chi$ in $U$. 
Similarily, if $T$ is a current of bidegree $(n+1,n+1)$ we define the
pushforward of $T$ by the same formula, but taking $\chi$ to be a function.
Since
$$
i\ddbar\pi_*(T)=\pi_*(i\ddbar T)
$$
it is clear
that $\pi_{*}(T)$ is subharmonic if $T$ is a subharmonic current of
bidegree $(n,n)$. 

If $T$ is an $(n,n)$-differential form with, say, bounded coefficients, the
pushforward of $T$ is a function whose value at a point $t$ equals
$$
\int_{\{t\}\times\Cn_z}T.
$$
Clearly, the pushforward only depends on the component of $T$ of
bidegree $(n,n)$ in $z$. Conversely, let $\kappa$ be a form of bidegree
$(n,n)$ in $z$, with coefficients depending on $t$. It follows from
the above that to prove that the  function
$$
\int_{\{t\}\times\Cn_z}\kappa
$$
is subharmonic it suffices to find a subharmonic form $T$ of
bidimension $(1,1)$ which  is compactly supported in the fiber
direction and  whose component of bidegree $(n,n)$ in $z$ equals
$\kappa$.

In order for this argument to work it is crucial that $T$ be globally
defined and compactly
supported in the fiber direction (or at least satisfies integrablility
conditions). The currents that we will encounter later are however
only defined in some pseudoconvex domain. To get globally defined
currents we  extend by 0 in the complement of the pseudoconvex
domain. This of course introduces a discontinuity which gives an extra
contribution to take into account when computing $i\ddbar T$ in the
sense of distributions. 
 The local calculations needed are summarized in
the following lemma, which is a variant of a by now standard method to
prove $L^2$-estimates for the $\dbar$-equation, see \cite{Horm} p 103 . 
\begin{lma} Let $\rho$ be a smooth real valued function in an open set
  $U$ in
  $\Cn$. Assume that $\partial\rho\neq 0$ on $S=\{z, \rho(z)=0 \}$, so
  that $S$ is a smooth real hypersurface. 
Let $T$ be a  real differential form of bidimension $(1,1)$ defined where
  $\rho<0$, with coefficients extending smoothly up to $S$. Assume 
$$
\partial\rho\wedge T
$$
vanishes on S, and that 
$$
\partial\rho\wedge\dbar\rho\wedge T
$$
vansihes to second order on $S$.
Extend
  $T$ to a current $\tilde T$ in $U$ by putting $\tilde T=0$
  where $\rho>0$. Then
\begin{equation}
i\ddbar\tilde T= \chi_{\rho<0}\quad i\ddbar T +\frac{i\ddbar\rho\wedge T
  dS}{|\partial\rho|} ,
\end{equation}
where $dS$ is surface measure on $S$ and $\chi$ is a characteristic function.
\end{lma}
In particular, even though it is not assumed that all of $T$, but only
certain components of $T$,  vanish on $S$, the contribution coming
from the discontinuity is a measure, and not, as might be expected, a
current of order 1. 
\begin{proof}
The hypotheses on $T$ mean that
\begin{equation}
\sum\rho_j T_{j\bar k}=\rho c_k,
\end{equation}
where $\sum c_k\rho_{\bar k}$ vanishes on $S$. Therefore, on $S$,
$$
0=\sum\frac{\partial}{\partial \bar z_k}(\rho c_k)=\sum\rho_j\frac{\partial
  T_{j\bar k}}{\partial\bar z_k} +\sum\rho_{j\bar k}T_{j\bar k}
$$
so
\begin{equation}
-\sum\rho_j\frac{\partial T_{j\bar k}}{\partial\bar z_k}=\sum
\rho_{j\bar k}T_{j\bar k}.
\end{equation}
Let $w$ be a smooth
function of compact support in $U$. Then, using the divergence  theorem
and  writing $T_{j \bar k}$ for
the components of $T$, we find that
$$
\int_{\rho<0}i\ddbar w\wedge T=\int_{\rho<0}\sum w_{j\bar k}T_{j\bar
  k}=
$$
$$
=\int_{\rho=0}\sum \rho_j w_{\bar k}T_{j\bar k}dS/|\partial\rho|
-\int_{\rho<0}\sum 
\frac{\partial w}{\partial z_j}\frac{\partial T_{j \bar
    k}}{\partial\bar z_k}.
$$
By equation (4.2) the boundary integral  vanishes. Applying the 
divergence theorem once more to the second integral we get 
 $$
\int_{\rho<0}w\sum\frac{\partial^2 T_{j\bar k}}{\partial
  z_j\partial\bar z_k}-\int_{\rho=0}w\sum\rho_j\frac{\partial T_{j\bar
    k}}{\partial\bar z_k}dS/|\partial\rho| .
$$
We then use (4.3) in the new boundary integral and find 
$$
\int_{\rho<0}i\ddbar w\wedge T=\int_{\rho<0}wi\ddbar T +\int_{\rho=0}\sum
\rho_{j\bar k}T_{j\bar k}dS/|\partial\rho|.
$$
This completes the proof of the lemma. 
\end{proof} 
 
\section{Second proof of Theorem 1.1}

Again, we first consider the situation described at the beginning of
section 2. 
As before, our starting point is the fact that the function
$$
u=\partial^{\phi}_t K_t
$$
is for fixed $t$ orthogonal to the space of holomorphic functions in
$A^2_t$. We now put
$$
k_t =K_t d\zeta_1\wedge...d\zeta_n,
$$
so that $k_t$ is $K_t$ interpreted as an $(n,0)$-form and, slightly abusively, define
$$
\partial^{\phi}_t k_t=\partial^{\phi}_t K_t\, d\zeta_1\wedge ...d\zeta_n.
$$

Since $\dbar$ has closed range, the orthogonal complement of the
kernel of $\dbar$ equals the range of $\dbar^*$. Therefore
$\partial^{\phi}_t k_t =\dbar^*\alpha$ for some 
 form $\alpha$ ( in $\zeta$) of bidegree $(n,1)$ which can also be taken to be
 $\dbar$-closed (and is then uniquely determined). By an
 argument similar to Lemma 2.1, 
 $\alpha$ depends smoothly on $t$. Write 
$\alpha=\sum\alpha_jd\bar \zeta_j\wedge d\zeta$.
Since $\alpha$ lies in the domain of $\dbar^*$, $\alpha$ satisfies the
$\dbar$-Neumann boundary condition $\sum \alpha_j\rho_j=0$ on the
boundary of $V$.

 Put $\gamma=\sum
\alpha_jd\hat{\zeta_j}$, where $d\hat\zeta_j$ stands for the wedge
product of all $d\zeta_k$:s except $d\zeta_j$, with a sign so that
$$
d\zeta_j\wedge d\hat\zeta_j=d\zeta_1\wedge...d\zeta_n.
$$
For later reference we note that the
$\dbar$-Neumann boundary condition on $\alpha$
translates to $\partial\rho\wedge\gamma=0$ on $\partial V$.
Put $g= dt\wedge\gamma +k_t$ and let $\partial^{\phi}=e^{\phi}\partial
e^{-\phi}$ be a twisted $\partial$-operator. 
The equation
$$
\partial^{\phi}_t k_t =\dbar^*\alpha
$$
is equivalent to
$$
\partial^{\phi}g=0. 
$$
We claim that the form $T$ defined as 
$$
T=c_n g\wedge\bar ge^{-\phi},
$$
 where $c_n$ is a constant of modulus 1 chosen so that $T$ is
 positive, for $\zeta$ in $V$ and $T=0$ for $\zeta$ outside of $V$ is a subharmonic
form. Since the component of $T$ of bidegree $(n,n)$ in $\zeta$ equals
$$
\kappa_t =c_n k_t\wedge\bar k_t
$$
it then  follows that 
$$
K_t(z,z)=\int \kappa_t
$$
is a subharmonic function of
$t$. 

To prove the subharmonicity of $T$ we first compute $i\ddbar T$ for
$\zeta$ inside of $V$. We use  the product rule
$$
\partial (a\wedge \bar b \,e^{-\phi})=\partial^{\phi} a\wedge \bar{b}\,
e^{-\phi} +(-1)^{\deg a}a\wedge\overline{\dbar b} \,e^{-\phi},
$$
and a similar rule for applying $\dbar$. Remembering that
$\partial^{\phi}g=0$ we get
\begin{equation}
i\ddbar T= c_n i\partial^{\phi}\dbar g\wedge\bar g \,e^{-\phi} + c_ni\dbar
g\wedge\overline{\dbar g} \,e^{-\phi}.
\end{equation}
From the commutation rule
$$
(\partial^{\phi}\dbar +\dbar\partial^{\phi})g=\ddbar\phi\wedge g,
$$
together with $\partial^{\phi} g=0$ it follows that the first term on
the right hand side can be written
$$
i\ddbar\phi\wedge T.
$$
This term is therefore nonnegative since $\phi$ is plurisubharmonic. 
To analyse the second term we introduce the notation $\zeta_0$ for $t$
and $\alpha_0$ for $-K_t$, so that $g$ can be written
$$
-\sum_0^n \alpha_jd\hat \zeta_j.
$$
The second term equals
$$
\sum_{j k}\frac{\partial\alpha_j}{\partial
  \bar\zeta_k}\overline{\frac{\partial\alpha_k}{\partial \bar\zeta_j}}.
$$
Here the indices run from 0 to $n$. Consider first the part of the sum
where both indices are greater than 0. Since the form $\alpha$ is
$\dbar$-closed for fixed $t$ this part equals
$$
\sum_1^n |\frac{\partial\alpha_j}{\partial
  \bar\zeta_k}|^2
$$
multiplied by the volume form $d\lambda$. Evidently, the part of the sum where
both indices are 0 equals
$$
|\frac{\partial\alpha_0}{\partial\bar\zeta_0}|^2 d\lambda.
$$
Finally, the terms in the sum when precisely one of the indices are 0
vanish since $\alpha_0=K_t$ is a holomorphic function of $\zeta$.

In conclusion, $i\ddbar T\geq 0$ for $\zeta$ in $V$. It now remains to
compute the contribution to $i\ddbar T$  which comes from cutting off
$T$ outside of $V$.
We  apply Lemma 4.1 to our current $T=c_n g\wedge\bar g e^{-\phi}$
and $\rho$ equal to the defining function of $V$. Then $\rho$ is
independent of $t=\zeta_0$, so 
$$\partial\rho\wedge
g=\partial\rho\wedge dt\wedge\gamma=0
$$
on $U\times\partial V$ since $\partial\rho\wedge\gamma=0$ on $\partial
V$.  
Hence the hypotheses of Lemma 4.1 are fulfilled. Since $V$ is
pseudoconvex it follows that $ic_n\ddbar\rho\wedge g\wedge\bar g$ is
non-negative on $\partial V$ so
$$
i\ddbar T\geq 0.
$$
In conclusion, $T$ is a subharmonic current and it follows that $K_t$
is a subharmonic function of $t$ for $z$ fixed. The rest of the proof
of Theorem 1.1 runs as before.

\section{Singularities of plurisubharmonic functions}

We first recall the definitions and basic properties of Lelong numbers
( our basic reference for these matters is \cite{Kis2}). If $\phi$ is
a plurisubharmonic function in an open set $U$ in $\Cn$ and $a$ is a
point in $U$, the {\it Lelong number} of $\phi$ at $a$ is
\begin{equation} 
\gamma(\phi,a)=\lim_{r\rightarrow 0}\,\,(\log r)^{-1}\sup_{|z-a|=r} \phi(z).
\end{equation}
Equivalently (see \cite{Kis2} p 176) we may introduce the mean value
of $\phi$ over the sphere centered at $a$ with radius $r$,
$M(\phi,a,r)$ and put
\begin{equation}
\gamma(\phi,a)=\lim_{r\rightarrow 0}\,\, (\log r)^{-1}M(\phi,a,r).
\end{equation}
The Lelong number measures the strength of the singularity of $\phi$
at $a$. If  $\gamma(\phi,a)>\tau$  then 
$$
\phi(z)\leq \tau \log|z-a|
$$
for $z$ close to $a$. 

In the one variable case we can decompose a subharmonic function
locally as a sum of  a harmonic part and a potential
$$
p(z) = \int\log |z-\zeta|d\mu(\zeta)
$$
where $\mu=1/(2\pi)\Delta \phi$. It is easy to verify that the Lelong
number is then equal to $\mu(\{a\})$. Using the potential it is also
easy to see that, in the one variable case, the
Lelong number at $a$  is greater than or equal to one if and only if
$e^{-2\phi}$ is not integrable over any neighbourhood of $a$.

In any dimension one defines  $\iota(\phi,a)$, the {\it integrability index} of $\phi$
at $a$,  as the infimum of all
positive numbers $t$ such that
$$
e^{-2\phi/t}
$$
is locally integrable in some neighbourhood of $a$. By a theorem of
Skoda (\cite{Skoda}), the inequality
$$
\iota(\phi,a)\leq \gamma(\phi,a)\leq n\iota(\phi,a)
$$
holds in any dimension. The left inequality here (which is the hard
part) says that if the Lelong number of $\phi$ at $a$ is strictly
smaller than 1, then $e^{-2\phi}$ is locally integrable near $a$.

Let $\Omega$ be a domain in $\Cn$ and let $\phi$ be a
plurisubharmonic function in $\Omega$. We consider the Bergman kernel
$K(z,z)$ for $A^2(\Omega,\phi)$. It is clear that if $a$ is a point
$\Omega$ and $e^{-\phi}$ is not integrable in any neighbourhood of
$a$, then any holomorphic function in  $A^2(\Omega,\phi)$ must vanish
at $a$, so in particular $K(a,a)=0$. Conversely, if $\Omega$ is bounded
and  $e^{-\phi}$
is integrable in some neighbourhood of $a$ then a standard application
of Hörmander's $L^2$-estimates shows that there exists some function
in $A^2(\Omega,\phi)$ which does not vanish at $a$. Since $K(a,a)$
equals the supremum of the modulus squared of all functions in
$A^2(\Omega,\phi)$ of norm 1, it follows that $K(a,a)>0$ in that
case. Thus, at least if $\Omega$ is bounded, the set where $\log K
=-\infty$ is precisely equal to the nonintegrability locus of
$e^{-\phi}$. 

For $z$ in $\Omega$ and $w$ in $\Cn$ we now consider the restriction
of $\phi$ to the complex line through $z$ determined by $w$ 
$$
\phi_{z,w}(\lambda)=\phi(z+\lambda w).
$$
For any fixed $z$ in $\Omega$ $\phi_{z,w}$ is defined for $\lambda$ in
the unit disk, if $w$ is small enough. Let $K_{z,w}(0,0)$ be the
Bergman kernel for the unit disk, with Lebesgue measure  normalized so that the
total area is one, equipped  with the weight function $2\phi_{z,w}$. By
the above, $K_{z,w}(0,0)=0$ if and only if the Lelong number of
$\phi_{z,w}$ at the origin is at least 1. By Theorem 1.1, $\log
K_{z,w}$ is a plurisubharmonic function, so for fixed $z$ the set of $w$ where  
it equals $-\infty$ is either pluripolar or contains a neighbourhood
of the origin. Thus, if the Lelong number at the origin of one single
slice function is smaller than 1, it must be smaller than 1 for all
slices outside a pluripolar set. 

It follows that the Lelong numbers of
all slices outside a pluripolar set are equal. This common value
also equals the Lelong number of $\phi$ at $z$.  To see this, first
note that by the first definition of Lelong number in terms of supremum
over spheres, it follows that the Lelong number for the restriction of
$\phi$ to any line through $z$ must be at least as big as the
$n$-dimensional Lelong number at $z$. The converse inequality follows
if we use the second definition of Lelong numbers in terms of mean
values over spheres, and apply Fatou's lemma. To avoid the
consideration of exceptional lines we now introduce the function
$$
\phi_\epsilon(z)=\frac{1}{2}\int_{|w|=\epsilon}\log K_{z,w}(0,0)dS(w),
$$
where the surface measure $dS$ is normalized so that the sphere has
total measure equal to 1. 
\begin{thm}
The function $\phi_\epsilon$ is well defined and plurisubharmonic in
the open set  $\Omega_\epsilon$ of points of $\Omega$ whose distance
to the 
boundary  is greater than $\epsilon$. The sequence $\phi_\epsilon$
decreases to $\phi$ as $\epsilon$ decreases to 0. The singularity set $S$
where $\phi_\epsilon=-\infty$ is for any $\epsilon>0$ equal to the
analytic set where the Lelong number of $\phi$ is at least 1. If the Lelong number
of $\phi$ at $z$ equals $\tau>1$, 
the Lelong number of $\phi_\epsilon$ at $z$ is at least equal to $\tau-1$. 
\end{thm}

\begin{proof}
Since $\log K_{z,w}$ is subharmonic with respect to $w$ t is clear
that $\phi_\epsilon$ decreases 
with $\epsilon$ to $\log K_{z,0}$. But $K_{z,0}$ is the Bergman kernel
at the origin 
for a normalized disk with a constant weight, $e^{-2\phi(z)}$, and so
equals $e^{2\phi(z)}$. Hence the limit of $\phi_\epsilon$ is equal to
$\phi$.  
If the Lelong number of $\phi$ at $z$ is smaller than 1 we have seen
above that $\log K_{z,w}(0,0)$ is not identically equal to $-\infty$
so its mean value over a sphere, $\phi_\epsilon$ is not equal to
$-\infty$. On the other hand we have also seen above that if
$\gamma(\phi,z)\geq 1$, then $\log K_{z,w}=-\infty$ for $w$ in a full
neighbourhod of 0, so $\phi_\epsilon(z)=-\infty$. Hence $S$ is equal
to the set where $\gamma(\phi,z)\geq 1$, which by Siu's analyticity
theorem, \cite{Siu}, is analytic. 

It remains only to prove the last statement of the theorem, so assume
0 lies in $\Omega$ and 
that $\gamma(\phi,0)=\tau>1$. Then, if $\tau'<\tau$,
$$
e^{-\phi(z)}\geq 1/|z|^{\tau'}
$$
if $|z|$ is small enough. For $w$ fixed and $h(\lambda)$ holomorphic we get
$$
\int_{|\lambda|<1}|h|^2 e^{-2\phi(z+\lambda w)}dm(\lambda)\geq\int_{|\lambda w|<
  |z|}|h|^2 e^{-2\phi(z+\lambda w)}dm(\lambda) \geq
$$
$$
\geq \int_{|\lambda|<1}|h|^2
/(|2z|^{2\tau'})dm(\lambda)\geq C |h(0)|^2/ |z|^{2(\tau'-1)}.
$$
Hence
$$
K_{z,w}(0,0)\leq C_1 |z|^{2(\tau'-1)}
$$
where the constant can be taken uniform for all $w$ of fixed modulus
equal to $\epsilon>0$. It follows that the Lelong number of
$\phi_\epsilon$ at $z$ is at least $\tau'-1$, and therefore at least
$\tau -1$ since $\tau'$ is an arbitrary number smaller than $\tau$.  
\end{proof}
The function $\phi_\epsilon$ thus ``attenuates the singularities'' of
$\phi$ in much the same way as Kiselman's construction in
\cite{Kis2}. ( Kiselman even gets that the Lelong number of the
constructed function {\it equals} $\tau -1$.)
In precisely the same way as in Kiselman, \cite{Kis3}, this
construction can be used to 
prove the Siu analyticity theorem. Let
$$
E_\tau=\{z;\gamma(\phi,z)\geq \tau\}.
$$
First, it follows from the Hörmander-Bombieri theorem that the
nonintegrability locus of any plurisubharmonic function is always
analytic. For a given plurisubharmonic function, $\phi$, and
$\delta>0$ we put, for some  choice of $\epsilon>0$
$$
\psi = 3n \phi_\epsilon/\delta.
$$
By Theorem 6.1, $\psi$ is finite at any point where
$\gamma(\phi,z)<1$, and therefore (see \cite{Horm}) $e^{-\psi}$ is
locally integrable near any such point. On the other hand $e^{\psi}$ is
not locally integrable near a point where
$\gamma(\phi,z)\geq(1+\delta)$ since the Lelong number of $\psi$ at
such a point is at least $3n$. Therefore we have, if $Z$ denotes the
nonintegrability locus of $e^{-\psi}$, that
$$
E_{1+\delta}\subset Z\subset E_1.
$$
Rescaling, we may of course for any $\tau>0$ and $\delta>0$ in a
similar way find an analytic set $Z_{\tau,\delta}$ such that
$$
E_\tau\subset Z_{\tau,\delta}\subset E_{\tau-\delta}.
$$
Hence $E_\tau$ equals the intersection of the analytic sets
$Z_{\tau,\delta}$ for $\delta>0$ and is therefore analytic. 

In a similar way we can consider, instead of restrictions of $\phi$ to
lines, the restriction of $\phi$ to $k$-dimensional subspaces. This
will give us a scale of ``Lelong numbers'' for $k=1, ...n$ that starts
with the classical Lelong number and ends with the integrability
index. 

We close this section by sketching an alternative way of relating
Lelong numbers to Bergman kernels, leading up to Theorem 1.3 of the
introduction. In \cite{Skoda} it is proved that if the Lelong number
of $\phi$ at $a$ is strictly smaller than 1, then $e^{-2\phi}$ is
locally integrable in some neighbourhood of $a$. Actually, Skoda's
proof of this fact gives a  bit more, namely that
$$
I(a) :=\int_{|z-a|<\delta} e^{-2\phi(z)}/|z-a|^{2n-2}dm(z)
$$
is also finite, if $\delta$ is small enough. ( The same argument as in
section 7 of \cite{Skoda} gives, with $d\sigma =\Delta\phi$  that
$$
\int_{|z|<r} e^{-2\phi(z)}/|z|^{2n-2}\leq C \int_{|z|<r,|x|<R}
|z|^{-2n-2}|z-x|^{-2n+\epsilon}d\sigma(x) dm(z),
$$
which is finite since 
$$
\int d\sigma(x)/|x|^{2n-2-\epsilon} 
$$
is finite.) On the other hand, $I(a)$ is comparable to the average of
$$
\int e^{-2\phi(a +\lambda w)} dm(\lambda)
$$
over all $w$ on a sphere, so $I(a)$ must be infinite if the Lelong
number of $\phi$ at $a$ is larger than or equal to 1. In conclusion
$$
\{a; I(a)=\infty\}=\{a, \gamma(\phi,a)\geq 1\}.
$$
We now introduce the plurisubharmonic function
$$
\psi(z,a) = \phi(z)+(n-1)\log|z-a|
$$
and let $K_a$ be the Bergman kernel for $\Omega$ with weight
$2\psi^a(z) =2\psi(z,a)$. It then follows that $\chi(a) =\log K_a(a,a)$
is plurisubharmonic and equal to $-\infty$ precisely where
$\gamma(\phi,\cdot)\geq 1$, so we have proved the first part of
Theorem 1.3 from the introduction. The last part of Theorem 1.3 follows from an
argument similar to the last part of the proof of Theorem
6.1. 
 
\section{ Plurisubharmonicity of potentials.} 

In this section we shall prove a  generalization of an earlier
result of Yamaguchi on the Robin function. Let $D$ be a smoothly
bounded pseudoconvex set in $\C_t^k\times\Cn_\zeta$ and let as before
$D_t$ be the $n$-dimensional slices of $D$. In this section we assume
$D$ has a smooth defining function $\rho(t,\zeta)$ such that
$\partial_{\zeta} \rho\neq 0$ on the boundary of $D$. In particular
all the fiber domains are smoothly bounded and have the same
topology. 

\begin{thm}
Let $K$ be a compact subset of $\Cn$ that is contained in $D_t$ for
all $t$ in an open set $U$. Let $\mu$ be a positive measure with support
in $K$. Let $u(t)$ be the negative of the energy of $\mu$ with respect to the Green
function $G_t$ of $D_t$ 
$$
u(t)=\int_{D_t} G_t(z,\zeta)d\mu(z)d\mu(\zeta).
$$
Then $u$ is plurisubharmonic in $U$. 
\end{thm}
Here we mean by the Green function the unique function vanishing on
the boundary and satisfying that
$\Delta_\zeta G $ is a unit point mass at $z$ .

The Green function $G$ of a domain $\Omega$ with pole at $z$ can be
written as the Newton 
kernel plus a smooth term
$$
G(z,\zeta)=-\frac{c_n}{|z-\zeta|^{2n-2}} + \psi(z,\zeta)
$$
where $\psi$ is harmonic in $z$ and in $\zeta$. The function
$$
\Lambda(z)=\psi(z,z)
$$
is called the Robin function of the domain $\Omega$. Let the  measure $\mu$
in Theorem 7.1 be a uniform mass distribution on a small ball centered
at the point $z$. $u$ is then equal to (the negative
of) the energy of
$\mu$ with respect to the Newton kernel plus the Robin function at $z$
of the domain $D_t$. Since
the Newtonian energy is independent of $t$ it follows that $\Lambda$
is plurisubharmonic as a function of $t$. Just like in the case of
 Bergman kernels this 
implies that even $\log\Lambda$ is subharmonic, if $n>1$. To see this, let $a$
be some complex number and consider the domain $D(a)$ with fibers
$$
D(a)_t =e^{at}D_t ,
$$
which, being a biholomorphic image if $D$, is still pseudoconvex. The
Robin function of $D(a)$ equals $e^{-(2n-2)\Re(at)}\Lambda$, so these
functions are subharmonic for any choice of $a$. It follows that
$\log\Lambda$ is subharmonic if $2n-2\neq 0$, i e if $n>1$. Finally,
we can again apply the Oka technique of variation of the domain (cf
the end of section 3) to conclude that if $\Lambda$ is the Robin
function of a fixed 
domain $\Omega$,  $\log\Lambda(z)$ is
plurisubharmonic as  a function of $z$ in $\Omega$. 

\begin{proof} We consider the Green function $G_t$ of
$D_t$ and let $g(t,z)=g_t(z)$ be the Green potential of $\mu$ in $D_t$. We
may assume that $\mu$ is given by a smooth density and it is then not
hard to see that $g$ is smooth up to the boundary in $D$. Let $\beta$
be the standard Euclidean Kähler form in $\Cn$ and set
$$
T=i\partial g\wedge\dbar g\wedge\beta_{n-1}
$$
in $D$ and $T=0$ outside of $D$.
(Here we use the notation $\omega_p=\omega^p/p!$ for $(1,1)$-forms
$\omega$.) Notice that $T$ is a nonnegative form and that $\pi_*(T)$
is given by
$$
\pi_*(T)(t)=\int_{D_t}|\partial g_t|^2=-\int_{D_t}\Delta g_t\, g_t=-u(t),
$$
where, as in Theorem 7.1, $u$ is the energy of $\mu$. 
Since $g$ vanishes on the boundary of $D$, $T$ satisfies the
hypotheses of Lemma 4.1. By Lemma 4.1
\begin{equation}
i\ddbar T\geq\chi_{D}i\ddbar T
\end{equation}
if D is pseudoconvex.  In $D$ we get 
\begin{equation}
i\ddbar T= -(i\ddbar g)^2\wedge\beta_{n-1}.
\end{equation}
Write
$$
i\partial\dbar g=i\partial_z\dbar_z g + i\partial_t\dbar_t g +
i\partial_t\dbar_z g 
+ i\partial_z\dbar_t g.
$$
Hence
$$
(i\ddbar g)^2\wedge\beta_{n-1}=
$$
$$
=2\Re (i\partial_t\dbar_tg\wedge i\partial_z\dbar_z g )\wedge\beta_{n-1} +
2\Re (i\partial_t\dbar_zg\wedge i\partial_z\dbar_t g) \wedge\beta_{n-1} =
$$
$$
= (2\Delta_t g\Delta_z g -2\sum|\frac{\partial^2g}{\partial z_j\partial
  \bar t}|^2)d\lambda.
$$
We thus find
$$
-i\ddbar u=\pi_*(i\ddbar T)\geq 
$$
$$
\geq2\left(\int_{D_t}-\Delta_t g\,\Delta_z g
+2\int_{D_t} \sum|\frac{\partial^2g}{\partial z_j\partial \bar t}|^2\right)
idt\wedge d\bar t \geq
$$
$$
\geq 2(-\Delta_t\int_{D_t}g_t d\mu)idt\wedge d\bar t =-2i\ddbar u.
$$
(Notice that we may move the Laplacian with respect to $t$ outside the
integral sign since $\mu$ is independent of $t$ and compactly
supported inside $D_t$.)
Thus $i\ddbar u\geq 0$, so $u$ is subharmonic and the proof of Theorem 7.1 
is complete. 
\end{proof}

Notice that the statement in Theorem 7.1 may be generalized to Green
functions for other elliptic equations, besides the Euclidean
Laplacian ( see also Yamaguchi and Levenberg \cite{Lev-Yam}). First, we may
replace the Euclidean metric by an arbitrary 
Kähler metric, with Kähler form $\omega$,  on $\Cn$, and consider the
Laplacian with respect to 
this metric. The same proof as above applies if we only replace the
Euclidean Kähler form $\beta$ by $\omega$. We may even go one step
further and consider elliptic operators of the form
$$
Lu =i\ddbar u\wedge \Omega
$$
where $\Omega$ is a closed positive form of bidegree $(n-1,n-1)$. 

It is also worth pointing out that the assumption on pseudoconvexity in
Theorem 7.1 can be relaxed. In the proof, convexity properties of the boundary of
$D$ only intervene in the application of Lemma 4.1, to conclude that
the form
$$
F = i\ddbar\rho\wedge i\partial g\wedge\dbar g\wedge\beta_{n-1} 
$$
is nonnegative on the boundary of $D$. Therefore we may replace the
hypothesis of pseudoconvexity in Theorem 7.1 by the hypothesis
``$F\geq 0$''. This is of course rather implicit, but to get an idea of how
 that condition relates  to pseudoconvexity we  can consider domains $D$ in
$\C\times\C^n$ of a special form. Let us assume e g that  
 the slices $D_t$ only depend on $\Re t$ and form an increasing family
 with respect to $\Re t$, 
 so that they are defined by inequalities
$$
D_t=\{z; v(z)<\Re t\}.
$$
When checking the positivity of the form $F$ above one may then replace
both $\rho$ and $g$ by $r=v-\Re t$, since  $\rho$ and $g$ are
positive multiples of $r$. We then see that, whereas the
pseudoconvexity of $D$ is equivalent to $v$ being {\it plurisubharmonic},
$F$ is positive if and 
only if $v$ is {\it subharmonic}. In particular this is a condition that
also makes sense  in $\R^n$.  In the next
section we shall briefly discuss analogs of the formalism of the last four
sections in $\R^n$. 

\section{Convexity properties of fiber integrals in $\Rn$}

We consider  $\R^{n+1}$ with the coordinates
 $(x_0, ... x_n)$ . When $\kappa$ is a function
with compact support or satisfying suitable integrability conditions,
we want to study convexity properties of the fiber integral
$$
\Phi(t)=\int_{x_0=t}\kappa dx_1...dx_n =\int_{x_0=t} \kappa.
$$
Just like in section 4 we shall arrange things so that 
$$
\kappa=T_{0 0}
$$
where $(T_{j k})$ is a matrix of functions. The basic fact of section
4, that the operation of pushforward of a form commutes with the
$i\ddbar$-operator, is now replaced by the following lemma.
\begin{lma}
Let $T=(T_{j k})$ be a matrix of $L^{\infty}$  functions in
$\R^{n+1}$. Suppose that for some $R>0$, $T$ vanishes when
$|(x_1,...x_n)|>R$.  Put
$$
\Phi(t)=\int_{x_0=t} T_{0 0}.
$$
If $T$ is smooth then
$$
\Phi''(t)=\int_{x_0=t}\sum_0^n \frac{\partial^2 T_{j k}}{\partial
  x_j\partial x_k}.
$$
If $T$ is not smooth the same formula holds in the sense of
distributions, if the right hand side is interpreted as the
distribution, $S$, whose action on a test function $\alpha$ is
$$
S.\alpha=\int_{\R^{n+1}}\sum_0^n \frac{\partial^2 T_{j k}}{\partial
  x_j\partial x_k}\alpha(x_0).
$$
\end{lma}
\begin{proof} If $T$ is smooth the first formula  is clear since the
  integral of any 
  term involving a 
  derivative with respect to a variable different from $x_0$
  vanishes. Hence
$$
\int_{x_0=t}\sum_0^n \frac{\partial^2 T_{j k}}{\partial
  x_j\partial x_k}=\int_{x_0=t} \frac{\partial^2 T_{0 0}}{\partial
  x_0\partial x_0}=\Phi''(t).
$$
The non-smooth case follows from the definition of distributional derivatives. 
\end{proof}
It is clear that the lemma holds even if $T$ does
not necessarily have 
compact support. It suffices that first and second order derivatives
of the coefficients of $T$ are integrable.
Later on we will also have use for a generalization of Lemma 4.1
\begin{lma}
Let $T=(T_{j k})$ be a matrix of   functions that are smooth up to the
boundary in a smoothly bounded domain $\Omega=\{\rho<0\}$ in
$\R^N$. Assume that
$$
\sum_j T_{j k}\rho_j=O(\rho)
$$
and
$$
\sum_j T_{j k}\rho_j\rho_k=O((\rho)^2)
$$
at the boundary of $D$. Extend the definition of $T$ to a matrix
$\tilde T$ in all of $\R^N$
by putting $\tilde T$ equal to 0 in the complement of $D$. Then we have in the
sense of distributions
$$
\sum \frac{\partial^2 \tilde T_{j k}}{\partial x_j\partial x_k}=
\chi_D\sum \frac{\partial^2  T_{j k}}{\partial x_j\partial x_k}+
\sum T_{j k}\rho_{j k}\frac{dS}{|d\rho|}.
$$
\end{lma}
\noindent The proof is essentially the same as the proof of lemma 4.1.
Let us now consider in particular matrices of the form
$$
T_{j k}=\gamma_j \gamma_k e^{-\phi}.
$$
To compute  derivatives of $T_{j k}$ we use the notation
$$
d_j = \partial/\partial x_j
$$
and
$$
d_j^\phi =e^\phi d_je^{-\phi}.
$$
We get
\begin{equation}
\frac{\partial^2 T_{j k}}{\partial  x_j\partial x_k}= \left(d_k\gamma_j
d_j\gamma_k + d_j^\phi d_k(\gamma_j)\gamma_k + d_j^\phi\gamma_j
d_k^\phi \gamma_k +\gamma_jd_jd_k^\phi \gamma_k\right)e^{-\phi}=
\end{equation}
$$
=\left(d_k\gamma_j
d_j\gamma_k + d_k d_j^\phi(\gamma_j)\gamma_k + d_j^\phi\gamma_j
d_k^\phi \gamma_k +\gamma_jd_jd_k^\phi \gamma_k +\phi_{j
  k}\gamma_j\gamma_k\right)e^{-\phi} 
$$
(where in the last line we have used the commutation relation
$$
d_j^\phi d_k=d_k d_j^\phi +\phi_{j k}.)
$$ 
It follows that if we assume 
$$
\sum d_k^\phi\gamma_k=0
$$
then
$$
\sum_0^n \frac{\partial^2 T_{j k}}{\partial x_j\partial
  x_k}=\left(\sum d_k\gamma_j d_j\gamma_k +|\sum d_k^\phi \gamma_k|^2
  +\phi_{j k}\gamma_j\gamma_k\right)e^{-\phi}.
$$
This identity can be used exactly as in section 5 to prove the real
Prekopa theorem. Let $\phi$ be a convex function and put
$$
k(t)=\gamma_0(t)=\left(\int_{x_0=t}e^{-\phi}\right)^{-1}.
$$
Since
$$
\int_{x_0=t} \gamma_0(x_0)e^{-\phi}=1
$$
it follows that
$$
\int_{x_0=t} d_0^\phi(\gamma_0) e^{-\phi}=0
$$
for any $t$. This implies that we can solve
$$
d_0^\phi \gamma_0=-\sum_1^n d_j^\phi \gamma_j e^{-\phi}
$$
and
$$
d_j\gamma_k=d_k\gamma_j
$$
with $\gamma$ and its first derivatives going rapidly to zero at
infinity (this is easiest to see for $n=1$, which is all that is
needed for the Prekopa theorem). If we then define $T_{j
  k}=\gamma_j\gamma_k e^{-\phi}$ as above it follows by Lemma 8.1 that
$$
k''(t)=\frac{d^2}{dt^2}\int_{x_0=t} k^2(x_0)e^{-\phi}= \int_{x_0=t}
\left(\sum d_k\gamma_j d_j \gamma_k +\phi_{j
    k}\gamma_j\gamma_k\right)e^{-\phi}.
$$
But, since $\gamma_0$ only depends on $x_0$, it follows just as in the
complex case that
$$
\sum d_k\gamma_j d_j \gamma_k= |d_0\gamma_0|^2
+\sum_1^n|d_k\gamma_j|^2\geq 0.
$$
Hence $k(t)$ is convex and it follows just like in section 3 that even
$\log c$ is convex, since replacing $\phi$ by $\phi +ax_0$ we see that
$k(t)e^{at}$ is convex for any choice of $a$. 

In the same way we can adapt the argument of section 7 to prove
convexity of Green potentials ( and hence the Robin function, see also
\cite{Cardaliaguet} who prove a stronger convexity property of the
Robin function), but in that case it is a little bit less evident
what the choice of $T$ should be. To explain this we shall first
discuss a general notion of subharmonic form in $\Rn$.

Consider the space, $F$, of differential forms on $\R^N_x\times\R^N_y$ whose
coefficients depend only on $x$. The usual exterior derivative, $d$,
preserves this space of forms. We introduce a new exterior derivative,
$d^{\#}$ on $F$ as
$$
 d^{\#}=\sum dy_j\wedge \partial/\partial x_j,
$$
where the partial derivative acts on the coefficients of a form (this
operator and the space $F$ are  not invariantly defined under changes of
coordinates). If we introduce the operator $\tau$ on $F$ by letting it
change $dx_j$ to $dy_j$ and vice versa, then $d^{\#}=\tau d\tau$ and it is
clear that $(d^{\#})^2=0$. We say that a form in $F$ is of bidegree
$(p,q)$ if its respective degrees in $dx$ and $dy$ are $p$ and
$q$. A $(p,p)$-form $\eta=\sum \eta_{I J} dx_I\wedge dy_J$ is symmetric if
$\eta_{I J}=\eta_{J I}$, or equivalently $\tau\eta=(-1)^p \eta$. Put
$$
\omega=\sum dx_j\wedge dy_j.
$$
A form  of bidegree $(N,N)$ is {\it positive} if it is a
nonnegative   multiple of $\omega_N:=\omega^N/N!$, and a general
symmetric form, $\eta$,
of bidegree 
$(p,p)$ is positive if
$$
a_1\wedge \tau a_1\wedge ...a_{N-p}\wedge\tau a_{N-p}\wedge\eta
$$
is positive for any choice of forms $a_j$ of bidegree $(1,0)$. It is
not hard to check that a form
$$
\sum a_{i j}dx_i\wedge dy_j
$$
is positive if and only if the matrix $(a_{i j})$ is positively
semidefinite. A smooth function $\phi$ is therefore convex precisely
when $dd^{\#}\phi$ is a positive form. It also follows that a positive
$(1,1)$-form can be written as a sum of forms of the type
$$
a\wedge\tau a,
$$
with $a$ of type $(1,0)$. Therefore the wedge product of a positive
form with a positive $(1,1)$-form is again positive. Similarily if we define
$dv_{j k}$ as the wedge product of all differentials
except $dx_j$ and $dy_k$, ordered so that $dx_j\wedge
dy_k\wedge dv_{j k} =\omega_N$, then
$$
\mu=\sum a_{j k} dv_{j k}
$$
is also positive exactly when $(a_{j k})$ is nonnegative as  a
matrix. A form $T=\sum T_{j k} dv_{j k}$ in $F$ of bidegree
$(N-1,N-1)$ is {\it subharmonic} if 
$$
dd^{\#} T=\sum_0^{n} \frac{\partial^2 T_{j k}}{\partial x_j\partial
  x_k}\omega_{N} 
$$ 
is positive. 

With these definitions it is clear that to apply Lemma 8.1 to prove
convexity of fiber integrals, we must look for subharmonic forms of
bidegree $(n,n)$ in $\R^{n+1}$. Let $D$ be a smoothly bounded 
domain in $\R^{n+1}$ defined by an inequality $D=\{\rho<0\}$ where the
gradient of $\rho$ does not vanish on the boundary of $D$ and let
$D_t$ be the $n$-dimensional 
slices. We say that such a domain satisfies condition $(C)$ if 
$$
d\rho\wedge d^{\#}\rho\wedge d d^{\#}\rho\wedge\omega'_{n-1}
$$
is  positive  for $x$ on the boundary of $D$. This condition is
clearly satisfied if $D$ is 
convex and it also holds if the fibers $D_t$ are of the form
$$
D_t=\{x'=(x_1,...x_n); v(x')<x_0\}
$$
where $v$ is subharmonic.
As in the case of Theorem 7.1 
we assume that the gradient of $\rho$ with
respect to $x'$ is never 0 for $x_0$ in an open set $U$ , so that all
the slices are smoothly 
bounded and have the same topology. Let $G_t$  be the Green function
of $D_t$. We then have.  
\begin{thm}
Assume that $D$ satisfies condition $(C)$. 
Let $K$ be a compact subset of $\Rn$ that is contained in $D_t$ for
all $t$ in  $U$. Let $\mu$ be a positive measure with support
in $K$. Let $u(t)$ be the negative of the energy of $\mu$ with respect to the Green
function $G_t$ of $D_t$ 
$$
u(t)=\int_{D_t} G_t(x,\xi)d\mu(x)d\mu(\xi).
$$
Then $u$ is convex in $U$. 
\end{thm}
The proof of Theorem 8.3 runs in much the same way as the proof of
Theorem 7.1. Let $g_t$ be the Green potential of $\mu$ in $D_t$ and
put
$$
g(x_0,x')=g_{x_0}(x').
$$
Let
$$
\omega'=\sum_1^n dx_j\wedge dy_j,
$$
and put
$$
T=dg\wedge d^{\#} g\wedge\omega'_{n-1}=\sum_0^n T_{j k}dv_{j k},
$$
for $x$ in $D$, and $T=0$ outside of $D$.
Then $T_{0 0}=|dg_{x_0}|^2$, so that
$$
\int_{x_0=t} T_{0 0}=-u(t).
$$
By Lemma 4.1 the contribution we get from the discontinuity at the
boundary of $D$  when we compute $dd^{\#} T$ equals
$$
dd^{\#}\rho\wedge T dS/|d\rho|.
$$
If $D$ satisfies condition $(C)$, this expression is nonnegative
(since $dg$ is a positive multiple of $d\rho$ at the boundary of $D$). By Lemma
8.2 we
therefore have (using $d\omega'=d^{\#}\omega'=0$) that
$$
\sum \frac{\partial^2 T_{j k}}{\partial x_j\partial x_k}\omega_{n+1}
=d d^{\#}T\geq-(dd^{\#}g)^2\wedge \omega'_{n-1}.
$$
Applying Lemma 8.1 we now get as in the complex case
$$
-u''(t)\geq \int_{x_0=t}\sum_1^n|\frac{\partial^2 g}{\partial
  x_j\partial x_0}|^2 -2u''(t),
$$
and it follows that $u$ is convex. 

Let us finally consider the implications of Theorem 8.3 for the Robin
function. Again as in section 7 we  take $\mu$ to be a positive measure
of total mass 1 which is given by a constant density on a small ball
centered at a fixed point $x$ that we assume to be contained in all
the fibers $D_t$, for $t$ in some open set $U$. The energy integral
$u(t)$ then equals  $\Lambda_t(x) -c$ where $\Lambda$ is the Robin
function for $D_t$ and $c$ is a constant. It follows that the Robin
function is a convex function of $t$ if $D$ satisfies condition
$(C)$. Moreover, $\Lambda$ is strictly convex at any point $t$ such
that the expression
$$
dd^{\#}\rho\wedge T
$$
is strictly positive at some point of the boundary of $D_t$. Consider
now the situation when all the fibers are translates of one fixed
domain $\Omega$ in $\R^n$
$$
D_t =\Omega +ta
$$
with $a$ a fixed direction in $\R^n$. Then $\rho(x_0,x)=r(x -x_0a)$ where
$r$ is a defining function for $\Omega$. It follows from the Hopf
lemma that $dg$ is a strictly positive multiple of $d\rho$ at the
boundary of $D$, so
$$
dd^{\#}\rho\wedge T
$$
is a strictly positive multiple of 
$$
\nu=dd^{\#}\rho\wedge d\rho\wedge d^{\#}\rho\wedge\omega'_{n-1}.
$$
To check the positivity of this $(n+1,n+1)$-form, we pull it back
under the map $(x_0,x,y_0,y)\rightarrow (x_0,x-x_0a,y_0,y-y_0a)$. It
is then not hard to see $\nu$ is positive for any choice of $a$ if
$\rho$ is convex and that moreover $\nu$ is strictly positive at any
point where the Hessian of $r$ restricted to the null space of $dr$ is
strictly positive. If $\Omega$ is smoothly bounded and convex there
will always be at least some such point at the boundary and the Robin
function is therefore strictly convex. 
We have
therefore proved a special case of a result from \cite{Cardaliaguet}:
\begin{thm} Let $\Omega$ be a  smoothly bounded convex domain in $\R^n$
  and let $\Lambda$ be the Robin function of $\Omega$. Then $\Lambda$
  is strictly convex.
\end{thm}
 In \cite{Cardaliaguet} a stronger convexity of the Robin function is
 proved (namely the {\it harmonic radius}, $\Lambda^{-1/(n-2)}$,  is
 strongly concave), but Theorem 8.4 is already sufficient to prove the
 unicity of 
 the {\it harmonic center} of $\Omega$, i e the point where $\Lambda$
 attains its minimum.

\bigskip

\def\listing#1#2#3{{\sc #1}:\ {\it #2}, \ #3.}


\begin{thebibliography}{9999}

\bibitem{Lieb}\listing  {H J Brascamp and E H  Lieb} {On
  extensions of the Brunn-Minkowski and Prékopa-Leindler theorems,
  including inequalities for log concave functions, and with an
  application to the diffusion equation.}  {J. Functional Analysis  22
  (1976), no. 4, 366--389.} 

\bibitem{BBN}\listing{K Ball, F Barthe and A Naor}{ Entropy jumps in
    the presence of a spectral gap} { Duke Math. J.  119  (2003),
    no. 1, 41--63.}


\bibitem{B1} \listing{B Berndtsson}{  Prekopa's theorem and Kiselman's
  minimum principle for plurisubharmonic functions.}{   Math. Ann.
  312  (1998),  no. 4, 785--792.}

\bibitem{Br-Bu}\listing{ J Bruna and J Burgu\'es}{ Holomorphic
    approximation and estimates for the $\overline\partial$-equation
    on strictly pseudoconvex nonsmooth domains}{  Duke Math. J.  55
    (1987),  no. 3, 539--596 }
\bibitem{Cardaliaguet}\listing{ Cardaliaguet, Pierre; Tahraoui,
    Rabah}{ On the strict concavity of the harmonic radius in
    dimension $N\ge3$.}{  J. Math. Pures Appl. (9)  81  (2002),
    no. 3, 223--240.} 

\bibitem{Cordero}\listing{D Cordero-Erausquin}{On Berndtsson's
    generalization of Prekopa's theorem} { Math Z  249 No 2 (2005)}

\bibitem{2Cordero}\listing{D Cordero-Erausquin}{Santaló's inequality on
    $\mathbb C\sp n$ by complex interpolation.}{
    C. R. Math. Acad. Sci. Paris  334  (2002),  no. 9, 767--772.} 

\bibitem{Demailly}\listing{J P Demailly}{Estimations $L\sp{2}$ pour
    l'opérateur $\bar \partial $ d'un fibré vectoriel holomorphe
    semi-positif au-dessus d'une variété kählérienne
    complète.}{  Ann. Sci. École Norm. Sup. (4)  15
    (1982)  no. 3, 457--511.} 
\bibitem{Horm}\listing{ L Hörmander}  {$L^2$-estimates and existence
  theorems for the $\dbar$- operator}{ Acta Math 113 (1965)}

\bibitem{Kis}\listing{ C O Kiselman}{The partial Legendre transformation for
  plurisubharmonic 
functions}{ Invent. Math. 49 (1978), no. 2, 137--148}

\bibitem{Kis3}\listing{ C O Kiselman}{Densit\'e des fonctions
    plurisousharmoniques}{Bull. Soc. Math. France. 107 (1979), no.3,
    pp295-304}

\bibitem{Kis2}\listing{ C O Kiselman}{ Attenuating the singularities
    of plurisubharmonic functions}{.  Ann. Polon. Math.  60  (1994),
    no. 2, 173--197}
\bibitem{Lev-Yam}\listing{ N Levenberg and H Yamaguchi}{ Robin
    functions for complex manifolds and applications}{Manuscript 2004}
\bibitem{Yam-Lev}\listing{ N Levenberg and H Yamaguchi}{ The metric
    induced by the Robin function.} {  Mem. Amer. Math. Soc.  92
    (1991),  no. 448} 
\bibitem{M-Y}\listing{ F Maitani and H Yamaguchi}{Variation of Bergman
    metrics on Riemann surfaces }{       Mathematische Annalen, Volume
    330, Number 3       477 - 489}

\bibitem{Prekopa}\listing{ A Prekopa} { On logarithmic concave
    measures and functions .}{ Acad Sci
Math (Szeged) 34 (1973) 335-343.}

\bibitem{Siu}\listing{ Y-T Siu}{ Analyticity of sets associated to
    Lelong numbers and the extension of closed positive currents.} {
Invent. Math. 27 (1974), 53--156.}

\bibitem{Skoda}\listing{H Skoda}{Sous-ensembles analytiques d'ordre
    fini ou infini dans $C\sp{n}$.}{  Bull. Soc. Math. France  100
    (1972), 353--408. } 
\bibitem{Yam}\listing{H Yamaguchi}{Variations of pseudoconvex domains
    over $C\sp n$}{
Michigan Math. J. 36 (1989), no. 3, 415--457.}
\end{thebibliography}
\end{document}